\documentclass[12pt]{amsart}
\usepackage{amscd,amssymb,verbatim}
\usepackage{amssymb,amsmath,amsthm,epsfig}

\newcommand{\N}{\mathbb{N}}

\newcommand{\R}{\mathbb{R}}
\newcommand{\C}{\mathbb{C}}

\theoremstyle{plain}
\newtheorem{Thm}{Theorem}[section]

\newtheorem{Lem}[Thm]{Lemma}

\newtheorem{Cor}[Thm]{Corollary}
\newtheorem{Rmk}[Thm]{Remark}

\newtheorem{Def}[Thm]{Definition}
\theoremstyle{remark}

\begin{document}

\title{A new method for constructing invariant subspaces}
\author{George Androulakis}

\date{}
\maketitle

\noindent {\bf Abstract:} The method of compatible sequences is introduced in order to 
produce non-trivial (closed) invariant
subspaces of (bounded linear) operators. Also a topological tool is used which is 
new in the search of invariant subspaces: the extraction of continuous selections of lower
semicontinuous set valued functions.
The advantage of this method over previously known methods is that
if an operator acts on a reflexive Banach space
then it has a non-trivial invariant subspace if and only if there exist compatible
sequences (their definition refers to a fixed operator).  
Using compatible sequences a result of Aronszajn-Smith is proved for reflexive Banach spaces.
Also it is shown that if $X$ be a  reflexive Banach space, $T \in {\mathcal L} (X)$, 
and $A$ is any closed ball of $X$,  
then either there exists $v \in A$ such that $Tv=0$, or there exists $v \in A$ such that
$\overline{\text{Span}}\, \text{Orb}_T (Tv)$
is a non-trivial invariant subspace of $T$, or 
$A \subseteq \overline{\text{Span}}\, \{ T^k x_{\ell} : \ell \in {\mathbb N}, 1 \leq k \leq \ell \} $
for every $(x_n)_n \in A^{\mathbb N}$.

\section{Motivation}

Let $X$ be a Banach space (always of dimension larger than $1$) and $T \in {\mathcal L}\, (X)$ (i.e.
$T:X \to X$ is  a (bounded linear) operator).
A (closed) subspace $Y$ of $X$ is called an invariant subspace of
$T$ if $T(Y) \subseteq Y$. Also $Y$ is called  non-trivial if
$\{ 0 \} \subsetneqq Y \subsetneqq X$.
The Invariant subspace problem asks whether every operator on a complex Hilbert space
has a non-trivial invariant subspace. This problem has its origins
approximately in 1935 when (according to \cite{AS}) J. von Neumann proved (unpublished) that
every
compact operator on a separable infinite dimensional complex Hilbert space has a non-trivial
invariant subspace,
(the proof uses the spectral theorem for normal operators, see \cite{vN}). 
Since then the invariant subspace
problem has motivated enormous literature in operator theory. The books \cite{B1}, \cite{P},
\cite{RR}, \cite{BFP}, the lecture notes \cite{A} and \cite{H2}, and the survey papers \cite{EL}
and \cite{AA} are centered around the invariant subspace problem. Related open problems and 
conjectures appear in \cite{AAST}. The invariant subspaces
appear in a natural way in prediction theory (see A.N. Kolmogorov \cite{K}, and N. Wiener
\cite{Wi}), and its interpretation in terms of unitary operators in a Hilbert space.

Let's recall some basic terminology and elementary facts about invariant subspaces: 
Let $X$ be a Banach space and $T \in {\mathcal L}\, (X)$.
The operator $T$ is called {\em non-transitive}
if it has a non-trivial invariant subspace. For $x \in X \backslash \{ 0 \}$ let
$\text{Orb}_T(x) = \{ T^n x : n \in \N \cup \{ 0 \} \}$ be the
orbit of $x$ under $T$. The closed linear space generated by $\text{Orb}_T(x)$,
$\overline{\text{Span}}\, \text{Orb}_T(x)$, is a invariant
subspace of $T$. The vector $x \in X \backslash \{ 0 \}$ is called {\em non-cyclic}
for $T$ if $\overline{ \text{Span}}\, \text{Orb}_T(x) \not = X$. If the operator $T$ is
non-transitive then it has a non-cyclic vector. Also if $X$ is a
non-separable Banach space and $T$ is an operator on $X$ then every non-zero vector
in $X$ is non-cyclic for $T$, thus $T$ is non-transitive.
If $X$ is a finite dimensional (always of dimension larger than $1$) complex
 or odd dimensional real Banach space then the Fundamental Theorem of Algebra
gives that $T$ has a non-zero eigenvector which spans a non-trivial invariant subspace.
If $X$ is an even dimensional real Banach space with dimension larger than $2$ then by considering
its complexification it can be proved that every operator
$T \in {\mathcal L}\, (X)$ has a non-trivial invariant subspace. On $\R^2$ the rotation by $\pi/2$
radians does not have non-trivial invariant subspaces.
If $T \not = 0$ and $Y:=\text{ker}\, (T) \not = \{ 0 \}$ then $Y$ is a non-trivial
(hyper-)invariant subspace of $T$ (i.e. $SY \subseteq Y$ for all
$S \in \{ T \} ' = \{ A \in {\mathcal L}\, (X): AT=TA \}$, the commutant of $T$).
If $T \not = 0$ and $Y:= \overline{\text{Ran}\, (T)} \not = X$
then $Y$ is a non-trivial hyper-invariant subspace of
$X$.

The key notion in this article is the notion of compatible sequences:

\begin{Def} \label{D:compatible}
Let $X$ be a Banach space and $T \in {\mathcal L}\, (X)$. Two bounded sequences 
$(x_n)_{n \in \N} \subset X$ and $(x_n^*)_{n \in \N} \subset X^*$ are called compatible
if 
\begin{itemize}
\item[(a)] No subsequence of $(x_n)_n$ converges weakly to zero.
\item[(b)] No subsequence of $(x_n^*)_n$ converges weakly to zero.
\item[(c)] There exists a sequence $(\varepsilon_n)_n$ of positive numbers converging to zero
such that 
\begin{equation} \label{mainequation}
 |x_m^* T^k x_{\ell}|< \varepsilon_m \text{ for all }1 \leq k \leq \ell \leq m.
\end{equation} 
\end{itemize}
\end{Def}

\begin{Rmk} \label{compatible}
Let $X$ be a reflexive Banach space and $T \in {\mathcal L}\, (X)$. $T$ has a non-trivial 
invariant subspace if and only if there exists a pair of compatible sequences.
\end{Rmk}

Indeed, let $(x_n)_n$ and $(x_n^*)_n$ be compatible sequences as in Definition~\ref{D:compatible}.
The reflexivity of $X$ implies 
the weak compactness of the closed balls, hence, since the sequences $(x_n)_n$ 
and $(x_n^*)_n$ are bounded, by the Theorem of Eberlein-\v{S}mulian,
pass to subsequences of $(x_n)_n$ and $(x_n^*)$ and
relabel to assume that there exists $v \in X$  and $v^* \in X^*$ 
such that
$(x_n)$ converges weakly to $v$ and $(x_n^*)_n$ converges weakly
to $v^*$. Notice that
(\ref{mainequation}) is still valid after the relabeling if one passes to the {\em same}
subsequence of $(x_n)_n$ and $(x_n^*)_n$ before relabeling. Now fix $1 \leq k \leq \ell$
and let $m \to \infty$ in (\ref{mainequation}) to obtain
\begin{equation} \label{2}
v^* T^k x_\ell =0 \text{ for all }1 \leq k \leq \ell.
\end{equation}
Then fix $k \in \N$ and let $\ell \to \infty$ in (\ref{2}) to obtain
\begin{equation} \label{3}
v^* T^k v =0 \text{ for all }k \in \N .
\end{equation}
By (a), $v \not = 0$. Thus, if $Tv=0$ then $\overline{\text{Span}}\, \text{Orb}_T(v)$ is a non-trivial
invariant subspace of $T$. If $Tv \not = 0$ then $\overline{\text{Span}}\, \text{Orb}_T(Tv)$
is a non-trivial invariant subspace of $T$ (it is non-trivial since by (\ref{3}) it is contained in the
kernel of $v^*$, $\text{ker}\, v^*$, and $\text{ker}\, v^* \not = X$ since $v^* \not = 0$ by (b)).

On the other hand, if an operator 
$T \in {\mathcal L}\, (X)$ has a non-trivial invariant subspace $Y$, then for every 
$x \in Y \backslash \{ 0 \}$, the closed linear space generated by $\text{Orb}_T(x)$ is a non-trivial
invariant subspace of $T$, thus there exists a non-zero functional $x^* \in X^*$ such that 
$\text{Orb}_T(x) \subseteq \text{ker}\, x^*$. Hence the sequences
$(x_n)_n$ and $(x_n^*)_n$ with $x_n=x$ and $x_n^*=x^*$ (for all $n$) are compatible. 
This finishes the proof of Remark~\ref{compatible}.

The method that is introduced in the present article, consists of showing how the absence of a 
pair of compatible 
sequences in some regions $A \subset X$ and $B \subset X^*$ (cf Claim~1 in the proof of 
Theorem~\ref{Th:AS}) implies (under the additional assumptions of Theorem~\ref{Th:AS}: 
the operator 
is compact and the reflexive Banach space is over the complex field) that the operator must 
have a non-zero eigenvector. The motivation for introducing the method of compatible sequences
comes from the equation in line 11 from the end of page 542 of \cite{AE}. In this equation a 
quantity of the form
$x_n^*T^kx_n$ is considered for some sequences $(x_n)_n \subset X$ and 
$(x_n^*)_n \subseteq X^*$ which satisfy (a) and (b) of Definition~\ref{D:compatible} and converge 
weakly to some non-zero elements $x\in X$ and $x^* \in X^*$ respectively. In order to conclude 
that $x_n^* T^k x_n$ converges to $x^*T^kx$ as $n \to \infty$, the assumption of the 
compactness of the operator $T$ is necessary. In order to avoid this necessity, information about
the behavior of $x_m^*T^kx_{\ell}$ is needed (i.e. when the indices of the sequences $(x_n)_n$
and $(x_n^*)_n$ do not necessarily match). 

A similarity of the present method to the one used in \cite{L1} is that the  
Brouwer-Schauder-Tychonov's fixed point
Theorem plays an important role in the proof. Nevertheless, the necessary 
compactness for the application of the Brouwer-Schauder-Tychonov's Theorem is not obtained 
by the compactness of the operator (as in \cite{L1}), but instead is obtained from the reflexivity
of the Banach space (hence convex closed bounded sets with the weak topology are compact). 
The reader may then
wonder why the  compactness of the operator is still assumed in Theorem~\ref{Th:AS}? 
The proof of Theorem~\ref{Th:AS} reveals that under the assumption of the compactness 
of the operator, the function $\Phi$ defined in (\ref{10}) is lower 
semicontinuous even if the range space $X$ is endowed with 
the norm topology (nevertheless notice that the lower semicontinuity of $\Phi$ if $X$ is 
endowed with the
weak topology, is sufficient for the proof to work!). Thus it is reasonable to ask whether 
modifications of the proof of Theorem~\ref{Th:AS} can remove or relax the assumption 
of the compactness of the  operator in the statement of Theorem~\ref{Th:AS}.

\section{Some history}

This section contains some (non-exhaustive) history. Many important directions
will not be discussed here. For some of these the reader is referred to \cite{A}, \cite{AA},
\cite{B2}, \cite{BFP}, \cite{EL}, \cite{H2}, \cite{P}, \cite{RR}. The history here is organized into
five subsections: {\bf (A)} the Theorem of Aronszajn-Smith and some extensions;
{\bf (B)} the Theorem of Lomonosov and some extensions; {\bf (C)} subnormal, hyponormal
operators and dual algebras; {\bf (D)} the method of extremal vectors and 
{\bf (E)} examples of transitive operators.

{\bf (A)} In 1954 N. Aronszajn and  K. Smith \cite{AS} proved that
if $X$ is an infinite dimensional complex Banach space and $T \in {\mathcal L}\, (X)$ is completely
continuous then $T$ has a non-trivial invariant subspace.
A non-linear map is used in the proof: $X$ is assumed without loss of generality to be
strictly convex and for a finite dimensional subspace $Y$ of $X$, the ``metric projection''
$P_Y:X \to Y$ (a non-necessarily linear map) is defined by
$P_Y(x)$ to be the unique $y \in Y$ which minimizes $\| x- y' \|$ for $y' \in Y $.

In 1966 A. Bernstein and  A. Robinson \cite{BR} proved that
if $H$ is a complex Hilbert space, $T \in {\mathcal L}\, (H)$ is a polynomially compact operator
(i.e. for some non-zero polynomial $p$, $p(T)$ is compact) then $T$ has a non-trivial invariant
subspace. The proof uses non-standard analysis as well as techniques similar to \cite{AS}.
In 1966 P.R.  Halmos gave a proof of the same result by a similar method but
avoiding the non-standard tools, \cite{H1}.

In 1968 W. Arveson and J. Feldman, \cite{AF} proved the following:
Let $H$ be a Hilbert space, and $T \in {\mathcal L}\, (H)$ satisfies
$\| TP_n - P_nTP_n \| \to 0$ for some sequence $(P_n)$ of orthogonal projection operators which
converges strongly to the identity operator (such operators are called quasitriangular;
the terminology is due to Halmos). Assume also that the norm closed algebra generated
by $T$ and $1$ contains a non-zero compact operator. Then $T$ has a non-trivial invariant
subspace. In 1973 C. Pearcy and N. Salinas \cite{PS} proved that
if $T$ is a quasitriangular operator on a Hilbert space and ${\mathcal R}\, (T)$
-the norm closure of the rational functions of $T$- contains a non-zero compact operator
then there exists a non-trivial subspace invariant under all operators in
${\mathcal R}\, (T)$. If $X$ is a Banach space and ${\mathcal A}$ is an algebra in
${\mathcal L}\, (X)$
then ${\mathcal A}$ is called non-transitive if there exists a non-trivial subspace which is
invariant under every element of ${\mathcal A}$. (If
${\mathcal A} = \{ T \}'$ then ${\mathcal A}$ is non-transitive if and only if
$T$ has a non-trivial hyperinvariant subspace.)

{\bf (B)} In 1973 V. Lomonosov \cite{L1} proved the following celebrated result:
Let $X$ be a complex Banach space and $T$ be an operator on $X$ which is not a multiple
of the identity and commutes with some non-zero compact operator.
Then $T$ has a non-trivial hyper-invariant subspace. For the proof of this result,
the new idea that was introduced was the Schauder fixed point Theorem:
If $\Phi$ is a norm-continuous function defined on a closed convex subset $C$ of a normed
space, and $\Phi (C)$ is contained in a norm-compact subset $K$ of $C$,
then $\Phi$ has a fixed point.
Lomonosov's result was extended to real Banach spaces by N.D.~Hooker in 1981, \cite{Ho}.
A special
case of Lomonosov's Theorem with a short proof (still using the Schauder fixed point Theorem)
was given by M. Hilden in 1977, \cite{H}.

A huge amount of literature has been produced towards extending Lomonosov's technique.
One of the strongest results in this direction was proved by Lomonosov in 1991 \cite{L2}:
Let $X$ be a Banach space, ${\mathcal A}$ be a proper subalgebra of ${\mathcal L}\, (X)$,
(i.e. ${\mathcal A} \not = {\mathcal L}\, (X)$) which is weakly closed. Then there exists
$x^{**} \in X^{**}\backslash \{ 0 \} $ and $x^* \in X^*\backslash \{ 0 \} $ such that
$|x^{**}T^*x^*| \leq ||| T ||| \text{  for all }T \in {\mathcal A},$
where $||| T |||$ denotes the essential norm of $T$ i.e. its distance
from the space of compact operators. A corollary of this result is the following:
Let $X$ be a Banach space and ${\mathcal A}$ be a weakly closed
proper subalgebra of ${\mathcal L}\, (X)$ such that
there exists a net
$\{ A_\alpha \} \subseteq {\mathcal A}$ and a non-zero operator $A \in {\mathcal A}$ such that
$A_{\alpha}^* \to A^*$ weakly and $||| A_{\alpha} ||| \to 0$.
Then $\{ T^*: T \in {\mathcal A} \}$ is non-transitive.
This corollary is a generalization of W.~Burnside's Theorem on matrix algebras: Every proper
algebra of matrices over an algebraically closed field has a non-trivial invariant subspace.
The shortest proof of Burnside's Theorem is given in \cite{LR}.

In 1996 A. Simoni\v{c} proved the following Hilbert space analogue of the 1991
Lomonosov's result \cite{S}:
Let $H$ be a complex Hilbert space, and ${\mathcal A}$ be a weakly closed proper  subalgebra
of ${\mathcal L}\, (H)$. Then there exist $x,y \in H \backslash \{ 0 \}$ such that
$ |\text{Re}\, <Tx, y>| \leq ||| \text{Re}\, T ||| <x,y>$ for all  $T \in {\mathcal A}.$
As a corollary he obtained that every essentially selfadjoint operator (i.e. $T-T^*$ is a compact
operator) on an infinite dimensional real Hilbert space has a non-trivial  invariant subspace.

{\bf (C)} Recall that a subnormal operator on a Hilbert space is the restriction of a normal operator 
to an invariant subspace. In 1978 S.W. Brown \cite{Br1} proved that
every subnormal operator has a non-trivial invariant subspace.
Functional Calculus was one of the main tools in the proof, which introduced the theory of
dual algebras. 
For the developments on this theory see the book \cite{BFP}. One of the main results 
produced by the  theory of dual algebras is the 1986 result of Brown, B. Chevreau and Pearcy
\cite{BCP}, that every contraction on the Hilbert space whose spectrum contains the unit circle 
has a non-trivial invariant subspace (also see \cite{B}).
Recall that an operator $T$ on a Hilbert space is called hyponormal if $TT^* \leq T^*T$.
Every subnormal operator is hyponormal.
In 1987 Brown proved \cite{Br2} that
every hyponormal operator $T$ has a non-trivial invariant subspace whenever
$C(\sigma(T)) \not = R(\sigma(T))$ where for a compact $K \subset \C$, $C(K)$
denotes the continuous functions on $K$ and $R(K)$ denotes the closure
(in the $C(K)$ norm) of the rational functions on $K$ with poles outside of $K$).
The  main tools of the proof are the theory of dual algebras and  a result of \cite{Pu} on 
properties of hyponormal operators.

{\bf (D)} The technique of ``extremal vectors'' was introduced in 1998 by A. Ansari and P. Enflo 
\cite{AE}. Let $X$ be a reflexive Banach space,  $T$
be an operator on $X$ with dense range,  $\varepsilon \in
(0,1)$, $x_0 \in X$, $\| x_0 \| =1$.  For every $n \in \N$ let
$y_n \in T^{-n} \{ x \in X: \| x_0 - x \| \le \varepsilon \} $
with $ \| y_n \| = \inf \{ \| y \| : y \in T^{-n} \{ x \in X:
\| x_0 - x \| \le \varepsilon \} \} .$ Then $(y_n)$ is called
a sequence of extremal vectors of $T$ with respect to
$\varepsilon$ and $x_0$. If $x$ is a weak limit point of $(T^ny_n)_n$ and $T$ satisfies
additional assumptions it turns out that $x$ is hyper-non-cyclic
for $T$ (i.e. $x$ is non-cyclic for all operators commuting with $T$).
Using the technique of minimal vectors, a special case of the 1973 result of Lomonosov
is proved in \cite{AE}: Let $K$ be a non-zero compact operator on a Hilbert space.  Then
$K$ has a non-trivial hyper-invariant subspace.

{\bf (E)} An example of a Banach space that admits an operator without any non-trivial
invariant subspace, (i.e. a transitive operator) was given by Enflo in the 70's, \cite{En1}, \cite{En2}.
The technique was simplified in \cite{B1}. Further examples were given by
C.~Read, \cite{R1}, \cite{R2}, \cite{R3}, \cite{R4}.

\section{The method of compatible sequences}\label{sec2}

Along with the method of compatible sequences, a new topological ingredient will be used 
in this article in order to obtain non-trivial invariant subspaces of operators: selection theorems 
for set 
valued functions. Excellent references for this topic are \cite{RS} and \cite{HP}.
Let ${\mathcal X}$ and
${\mathcal Y}$ be topological spaces, ${\mathcal P}\, ({\mathcal Y})$ denote the power set of
${\mathcal Y}$ and $\phi: {\mathcal X} \to {\mathcal P}\, ({\mathcal Y})$ be a set valued function. 
If $\phi (x) \not = \emptyset$ for all $x \in {\mathcal X}$ then a function
$f:{\mathcal X} \to {\mathcal Y}$ is called a selection of $\phi$ if $f(x) \in \phi (x)$ for all
$x \in {\mathcal X}$. The set valued function $\phi$ is called lower semicontinuous (l.s.c.)
if for any open subset $V$ of ${\mathcal Y}$,
the set $\{ x \in {\mathcal X}: \phi (x) \cap V \not = \emptyset \}$ is open in ${\mathcal X}$.
In terms of convergence of nets, this definition is equivalent to: for all $x \in {\mathcal X}$,
for all $y \in \phi (x)$ and for all nets 
$(x_{\lambda})_{\lambda \in \Lambda} \subseteq {\mathcal X}$ with $x_{\lambda} \to x$,
obtain that there exists $y_{\lambda} \in \phi (x_{\lambda})$ such that $y_{\lambda} \to y$.
For the set valued map $\phi$ denote by $\overline{\phi}$ and $\overline{\text{conv}}\, \phi$ 
(if the values of $\phi$ are subsets of a topological vector space) the maps defined by
$$
\overline{\phi}(x) = \text{ the closure of the set } \phi(x); \phantom{aa}
\overline{\text{conv}}\, \phi (x)= \text{ the closure of the convex hull of } \phi (x).
$$
Recall that a Hausdorff topological space ${\mathcal X}$ is called paracompact if every open
covering of ${\mathcal X}$ has a locally finite refinement. 
Michael's selection Theorem, \cite[Theorem~$3.2''$]{M}, states that a Hausdorff topological space
${\mathcal X}$ is paracompact if and only if for any Banach space ${\mathcal Y}$, every
l.s.c. set valued function $\phi : {\mathcal X} \to {\mathcal P}\, ({\mathcal Y})$ whose
values are non-empty closed convex subsets of ${\mathcal Y}$ has a continuous selection. 
It should be mentioned that if one is interested in sufficient conditions for the existence 
 of continuous selections of
set valued functions, then the assumptions of the above result are not optimal. The assumption
of lower semicontinuity can be weakened to ``weak lower semicontinuity'' as it was observed 
by K. Przeslawski and L. Rybinski (see \cite{PR} for details). The assumption that ${\mathcal Y}$ 
is a Banach space can also be weakened. In fact in \cite{M} it was observed that the proof gives
that ${\mathcal Y}$ can be assumed to be an $F$ space (i.e. a complete metrizable topological
vector space). The metrizability of ${\mathcal Y}$ is not necessary, since the result was further 
improved in \cite{M3} as follows:

\begin{Thm} \label{selection}
Let ${\mathcal X}$ be a paracompact topological space, ${\mathcal Y}$ be a locally convex 
topological vector
space with the property that the closed convex hull of any compact set is compact, and 
$\phi: {\mathcal X} \to {\mathcal P}\, ({\mathcal Y})$ be an l.s.c. map such that 
$\cup \{ \phi(x): x \in {\mathcal X} \} $ is metrizable and $\phi(x)$ is a complete set for all 
$x \in {\mathcal X}$. Then the map $\overline{\text{conv}}\, \phi$ admits a continuous selection.  
\end{Thm}

Now the technique of compatible sequences is presented. It is applied in order to obtain the 
following result which is a special case of the result of \cite{AS}.

\begin{Thm} \label{Th:AS}
Every compact operator on a complex reflexive Banach space has a non-trivial invariant subspace.
\end{Thm}

\begin{proof}
Let $X$ be a complex reflexive Banach space and $T$ be a compact operator on $X$. 
It will be shown that $T$ has a non-trivial invariant subspace.
Without loss of generality assume that $X$ is separable (else the closed linear
span of the the orbit of any non-zero vector under $T$ is a non-trivial invariant
subspace of $T$). 

\medskip
\noindent
{\bf Claim 1:} There exist sets $A \subseteq X$, $B \subseteq X^*$ and sequences
$(x_n)_{n \in \N} \subseteq A$ and $(x_n^*)_{n \in \N} \subseteq B$ satisfying
\begin{itemize}
\item[(a)] $A$ is bounded and the origin does not belong to the weak closure of $A$.
\item[(b)] $B$ is bounded and the origin does not belong to the weak closure of $B$.
\item[(c)] For all integers $1 \leq k \leq \ell \leq m$,
\begin{equation} \label{1}
|x_m^* T^{k}x_{\ell} | < \frac{1}{m}.
\end{equation}
\end{itemize}

Once Claim~$1$ is established, the sequences $(x_n)_n$ and $(x_n^*)_n$ are compatible, thus 
by Remark~\ref{compatible}, $T$ has a non-trivial invariant subspace. This finishes the proof.
It remains to prove Claim~$1$.

The sets $A$ and $B$ that appear in (a) and (b) of Claim~$1$ are constructed as
follows: The set $A$ is any non-empty closed convex bounded subset of $X$ which does not 
contain the  origin and has a non-empty norm interior. 
Let $x_0$  be a point in the norm-interior of $A$. Let 
$$
B = \{ x^* \in X^* : \text{Re}\, x^*x \geq 0 \text{ for all }x \in A 
\text{ and }\text{Re}\, x^* x_0 =1 \} . 
$$
Since $A$ is convex, has non-empty norm-interior and does not contain the origin, by the 
Eidelheit separation theorem, $B$ is non-empty.
Obviously $B$ is weakly closed. In order to see that $B$ is bounded,
fix $x^* \in B$ and notice that since $x_0$ belongs to the norm-interior of $A$,
there exists $\varepsilon >0$ such that $x_0 -x \in A$ for all $x \in X$ with 
$\| x \| \leq \varepsilon$. Hence
$$
0 \leq \text{Re}\, x^* (x_0 -x) = \text{Re}\, x^* x_0 - \text{Re}\, x^* x=
 1- \text{Re}\, x^* x.
$$
Thus $\text{Re}\, x^* x \leq 1$ for all $x \in X$ with $\| x \| \leq \varepsilon$, 
which proves that $B$ is bounded. Hence $B$ is weakly compact, and 
the sets $A$, $B$ constructed above satisfy (a) and (b) of Claim~$1$.
It remains to construct sequences $(x_n)_n \subseteq A$ and $(x_n^*)_n \subseteq B$
satisfying (c) of Claim~$1$.

\medskip
\noindent
{\bf Claim 2:} Either there exists a sequence $(x_n)_{n \in \N} \subseteq A$ such that
for every $m \in \N$,
$$
B \cap \bigcap_{\ell =1}^m \bigcap_{k=1}^{\ell}
\{ x^* \in X^*: | x^*T^k x_{\ell}   | < \frac{1}{m} \}
\not = \emptyset,
$$
or $T$ has an eigenvalue.

Of course, if the first alternative of Claim~$2$ is valid then pick a sequence $(x_n^*)_{n \in \N}$
in $X^*$ with
$$
x_m^* \in B \cap \bigcap_{\ell =1}^m \bigcap_{k=1}^{\ell}
\{ x^* \in X^* : | x^* T^k x_{\ell} | < \frac{1}{m} \}
$$
and obviously $(x_n)_{n \in \N}$ and $(x_n^*)_{n \in \N }$ satisfy (c) of Claim~$1$ and
the proof finishes.
Also, if the second alternative of Claim~$2$ is valid then the proof finishes as well.
Thus it only remains to establish Claim~$2$.

Assume that the first alternative of Claim~$2$ is false. It will be shown that the second
alternative of Claim~$2$ is true. Since the first alternative of Claim~$2$ is false, for every
sequence $(x_n)_{n \in \N} \subseteq A$ there exists $m \in \N$ such that
$$
B \cap \bigcap_{\ell =1}^m \bigcap_{k=1}^{\ell}
\{ x^* \in X^* : |x^* T^k x_{\ell}  | < \frac{1}{m} \} = \emptyset ,
$$
or equivalently,
\begin{equation} \label{4}
B \subseteq \bigcup_{\ell =1}^m \bigcup_{k=1}^{\ell}
\{ x^* \in X^* : |x^* T^k x_{\ell} | \geq \frac{1}{m} \} .
\end{equation}
For $(x_n)_n \subseteq A$ and $m \in \N$ define
\begin{equation} \label{8}
Y((x_n)_{n=1}^m):= \text{Span}\, \{  T^k x_{\ell} : \text{ for } \ell \in \{ 1 , \ldots , m \}  \text{ and }
k \in \{ 1, \ldots , \ell \} \} .
\end{equation}

\medskip
\noindent
{\bf Claim 3:} Let any sequence $(x_n)_{n \in \N} \subseteq A$ and pick $m \in \N$ for which
(\ref{4}) is valid.
Then  $Y((x_n)_{n=1}^m) \cap A^o \not = \emptyset$, where $A^o$ denotes the norm-interior
of $A$. 

Indeed, let $(x_n)_{ n \in \N} \subseteq A$ 
and $m \in \N$ for which (\ref{4}) is valid. If $Y((x_n)_{n=1}^m)\cap A^o = \emptyset$, 
then since $Y((x_n)_{n=1}^m)$ and $A^o$ are 
convex sets and $A^o$ is non-empty open, by the Eidelheit separation theorem
there exists $x^* \in X^*$ such that the restriction of $x^*$ on 
$Y((x_n)_{n=1}^m)$ is equal to zero, (since $Y((x_n)_{n=1}^m)$ is a linear space), 
$\text{Re}\, x^* x \geq 0$ for all $x \in A$, and $\text{Re}\, x^* x > 0$ for all 
$x \in A^o$. In particular, $\text{Re}\, x^* x_0 > 0$.
Let $y^* = x^*/\text{Re}\, x^* x_0$. Then $y^* \in B$.
Hence (\ref{4}) implies that there exists $\ell \in \{ 1, \ldots , m \}$ and 
$k \in \{ 1, \ldots \ell \}$ such that $| y^* T^k x_{\ell}  | \geq \frac{1}{m}$. 
In particular, $y^* T^k x_{\ell} \not = 0$ which contradicts the fact that 
the restriction of $y^*$ on  $Y((x_n)_{n=1}^m)$ is equal to zero,
and finishes the proof of Claim~$3$.

Define a function $m: A^{\N} \to \N$ as follows: For $(x_n)_{n \in \N} \in A^{\N}$ let
\begin{equation} \label{9}
m((x_n)_{n \in \N}) = \min \{ m \in \N :
B \subseteq \bigcup_{\ell =1}^m \bigcup_{k=1}^{\ell}
\{ x^* \in X^* : |x^* T^k x_{\ell}  | \geq \frac{1}{m} \} \} .
\end{equation}
Define a set valued function $\Phi : A^{\N} \to {\mathcal P}\, (X)$ by 
\begin{equation} \label{10}
\Phi((x_n)_{n \in \N}) = Y((x_n)_{n=1}^{\mu}) \cap A^o
\text{ where } \mu =m((x_n)_{n \in \N}) .
\end{equation}
In the proof of Claim~$2$ the following has been achieved so far:

{\em If the first alternative of Claim~$2$ is false then the set valued function $\Phi$ defined
in (\ref{10}) has the property that $\Phi( (x_n)_n)$ is a non-empty convex set 
for all $(x_n)_n \in A^{\N}$; (the fact that for every $(x_n)_n \in A^{\N}$, $\Phi( (x_n)_n)$
is a non-empty set, is exactly the statement of Claim~$3$; the fact that $\Phi( (x_n)_n)$ is a 
convex set is obvious).}

For the remaining of the proof, endow $X$ (and thus $A$) with the weak topology and 
$X^{\N}$ (and thus $A^{\N}$) with the product topology of the weak topology. 
Then, since $X$ is reflexive, $A$ and $A^{\N}$ are compact topological spaces.

\medskip
\noindent
{\bf Claim 4:} The function $\Phi$ defined in (\ref{10}) is l.s.c.

Assume for the moment that Claim~$4$ has been shown. Now observe that the assumptions 
of Theorem~\ref{selection} are satisfied for $\Phi$: Indeed
$A^{\N}$ is a compact Hausdorff space (hence it is paracompact). 
Since $X$ is separable reflexive and $A$ is bounded, the weak topology on $A$ is metrizable.
Thus by Theorem~\ref{selection} 
there exists a continuous selection $F:A^{\N} \to X$ of $\overline{\Phi}$. 
For every $(x_n)_n \in A^{\N}$, if $\mu = m((x_n)_n)$ then
$$
\overline{\Phi}((x_n)_{n \in \N}) \subseteq A \cap  Y((x_n)_{n=1}^{\mu}) .
$$
Hence $F((x_n)_n) \in A$ for all $(x_n)_n \in A^{\N}$.  Define
$$
G: A^{\N} \to A^{\N} \text{ by } G((x_n)_n) =
(F((x_n)_n), F((x_n)_n) , \ldots ).
$$
Obviously $G$ is continuous since $F$ is continuous. Since $A^{\N}$ is a convex
compact subset of the topological vector space $X^{\N}$ (endowed with the product topology
of the weak topology of $X$), by the Brouwer-Schauder-Tychonov's Theorem,
(see for instance \cite[Corollary $16.52$]{AB}), $G$ has a
fixed point. If $(x_n)_n \in A^{\N}$ is a fixed point of $G$ then
$F((x_n)_{n \in \N})=x_1=x_2= \cdots$. Thus there exists $x \in A$ such that
$F((x_n)_n)=x$ where $x_n=x$ for all $n \in \N$. If $\mu =m((x_n)_n)$
then by the definition of $Y((x_n)_{n=1}^\mu)$, $F((x_n)_n)$ is a finite linear
combination of the vectors $\{ T^k x: k=1, \ldots , \mu \}$. Thus there exists a non-constant
polynomial $q$ such that $F((x_n)_n)= q(T)x$. Hence $q(T)x =x$, i.e.
$1 \in \sigma_p(q(T))$. Since $X$ is a complex Banach space, by the Fundamental Theorem of
Algebra,
$\sigma_p(q(T))=q(\sigma_p(T))$ hence $\sigma_p(T) \not = \emptyset$, thus $T$ has a
non-zero eigenvector which spans a non-trivial invariant subspace of $T$. Hence the second
alternative of Claim~$2$ was shown to be true. It only remains to establish Claim~$4$.

The proof of Claim~$4$ is based on the fact that the function 
$m: A^{\N} \to \N$ is an l.s.c. function. In order to see that $m$ is an l.s.c. function, let 
$((x_{j,n})_{n \in \N})_{j \in \N} \subseteq A^{\N}$ and $(x_n)_{n \in \N} \in A^{\N}$
such that $(x_{j,n})_{n \in \N} \to (x_n)_{n \in \N}$ as $j \to \infty$, i.e.
\begin{equation} \label{15}
x_{j,\ell} \to x_\ell \text{ weakly as }j \to \infty , \text{ for all }\ell \in \N .
\end{equation}
Let $\mu :=m((x_n^*)_{n \in \N})$ and for $j \in \N$, let $\mu_j := m((x_{j,n}^*)_{n \in \N})$ as
defined by (\ref{9}).
It is claimed that $\mu_j \geq \mu$ for all $j$ large enough (which will finish the proof that
the function $m$ is an l.s.c. function). Indeed, if $\mu_j \leq \mu -1$ for
infinitely many $j$'s, then by passing to a subsequence and relabeling assume that
there exists a constant $\mu_0 \leq \mu -1$ such that $\mu_j= \mu_0$ for all $j \in \N$. 
Thus by (\ref{9}),
\begin{equation} \label{13}
B  \subseteq \bigcup_{\ell =1}^{\mu_0} \bigcup_{k=1}^{\ell} \{ x^* \in X^*:
| x^* T^k x_{j,\ell}  | \geq \frac{1}{\mu_0} \} \text{ for all }j \in \N ,
\end{equation}
 yet
\begin{equation}\label{14}
B \not \subseteq \bigcup_{\ell =1}^{\mu_0} \bigcup_{k=1}^{\ell}
\{ x^* \in X^*: | x^* T^k x_{\ell} | \geq \frac{1}{\mu_0} \} .
\end{equation}
By (\ref{14}) there exists $b \in B$ such that
$$
| b T^k x_\ell | < \frac{1}{\mu_0} \text{ for all }\ell =1, \ldots , \mu_0 \text{ and }
k=1, \ldots , \ell .
$$
Then by (\ref{15}) there exists $j$ large enough such that
\begin{equation} \label{16}
|b T^k x_{j,\ell} | < \frac{1}{\mu_0} \text{ for all }\ell =1, \ldots , \mu_0 \text{ and }
k=1, \ldots , \ell .
\end{equation}
But (\ref{16}) contradicts (\ref{13}). The contradiction proves that  $\mu _j \geq \mu$ for all $j$
large enough, and by passing to a subsequence and relabeling, assume that $\mu_j \geq \mu$
for all $j \in \N$.

The proof of Claim~$4$ continues as follows: 
Let $((x_{j,n})_{n \in \N})_{j \in \N} \subseteq A^{\N}$ and 
$(x_n)_{n \in \N} \in A^{\N}$ satisfying (\ref{15}), $\mu :=m((x_n)_{n \in \N})$, 
$\mu _j := m((x_{j,n})_{n \in \N})$ (for $j \in \N$), and $y \in \Phi ((x_n)_{n \in \N})$.
By (\ref{8}) and (\ref{10}) there exist $(a_{k,\ell})_{\ell=1, k=1}^{m, \ell} \subset \C$ such that
$$
y = \sum_{\ell =1}^{\mu} \sum_{k=1}^\ell a_{k,\ell}T^k x_\ell .
$$
For $j \in \N$ let $Y_j:= Y((x_{j,n})_{n=1}^{\mu_j})$ as defined in (\ref{8}), and
$$
y_j:= \sum_{\ell=1}^{\mu} \sum_{k=1}^{\ell} a_{k,\ell}T^k x_{j,\ell} .
$$
Then, for all $j$'s large enough,  $\mu_j \geq \mu$ (since the function $m$ is an l.s.c. function),
and thus $y_j \in Y_j$. By (\ref{15}) and the compactness of $T$ obtain that   $y_j \to y$ in norm.
Thus $y_j \in A^o$ (since $y \in A^o$) and hence $y_j \in A^o$ and 
$y_j \in \Phi ((x_{j,n})_{n \in \N})$ for all 
$j$'s large enough.
Since $y_j \to y$ weakly as $j \to \infty$, the set valued function $\Phi$ is l.s.c.
which finishes the proof of Claim~$4$.
 \end{proof}

Notice that in the proof of Theorem~\ref{Th:AS}, the compactness of the operator $T$ was only 
used in the third line before the end of proof. Thus the compactness of $T$ was used to show that 
the map $\Phi$ defined in (\ref{10}) is l.s.c. In fact in the third line before the end of the proof 
of Theorem~\ref{Th:AS}, it is obtained that ``$y_j \to y$ {\em in norm}'' which implies that the 
function $\Phi$ is l.s.c. {\em even if the range space is equipped with
the norm topology} (while still the domain is equipped with the product of the weak topology).
Note that for the proof of Theorem~\ref{Th:AS} to work, it is only required that $\Phi$ is l.s.c.
when the range space is equipped with the weak topology.
Thus it is reasonable to ask whether the assumption of the compactness of the operator in the
above proof can be omitted, or replaced by a weaker assumption.

In the proof of Theorem~\ref{Th:AS}, under the assumption of the compactness of the operator 
$T$, the use of the selection Theorem~\ref{selection} was not necessary. Alternatively $A$ can be 
taken to be a closed strictly convex ball of $X$ which does not contain the origin
(without loss of generality $X$ can be renormed to be strictly convex) and it can be
shown that the function $F:A^{\N} \to A$ defined by $F((x_n)_n)$ to be the unique point of 
$A \cap Y((x_n)_{n=1}^{m((x_n)_n)})$ of the smallest norm, is a continuous selection of 
$\overline{\Phi}$. This choice of a continuous selection is motivated by \cite{AE} and \cite{AS}.

Another application of the method of compatible sequences is the following:

\begin{Thm} \label{Th:alternatives}
Let $X$ be a  reflexive Banach space, $T \in {\mathcal L} (X)$, and $A$ be any 
closed ball of $X$.  
Then either there exists $v \in A$ such that $Tv=0$, or there exists $v \in A$ such that
$\overline{\text{Span}}\, \text{Orb}_T (Tv)$
is a non-trivial invariant subspace of $T$, or 
$A \subseteq \overline{\text{Span}}\, \{ T^k x_\ell : \ell \in \N, 1 \leq k \leq \ell \} $
for every $(x_n)_{n \in \N} \in A^{\N}$.
\end{Thm}

Notice that if $X$ is a Banach space and $T \in {\mathcal L} (X)$ which does not have any 
non-trivial invariant subspace
then for any $x \in X\backslash \{ 0 \}$, the span of the set $\{ T^k x: k \in \N \}$ is dense in 
$X$. Thus Theorem~\ref{Th:alternatives} is a localized version of the following obvious 
observation: Every $T \in {\mathcal L} (X)$ either has a non-trivial 
invariant subspace or for every 
$x \in X \backslash \{ 0 \}$ the span of the set $\{ T^k x : k \in \N \}$ is dense in $X$.

\begin{proof}[Proof of Theorem~\ref{Th:alternatives}]
 Let $X$ be a reflexive  Banach space and $T$ be an operator on $X$. Let $A$ be any 
closed ball in $X$. If $A$ contains the origin then the first alternative of 
Theorem~\ref{Th:alternatives} is valid, thus assume that $A$ does not contain the origin.  
Let $A^o$ denote the norm-interior of $A$. 
Fix a weakly open neighborhood  $V$ of the origin of $X$. Observe that since 
$\{ y+V: y \in A^o \}$ is a weakly open cover of the weakly compact set $A$
(by the reflexivity of $X$, $A$ is weakly compact; note also that $A^o$ is dense in $A$), 
there exists a 
positive integer $n(V)$ and a finite set $\{ a(V,1), a(V,2), \ldots , a(V,n(V)) \} \subseteq A^o$
such that $\cup \{ a(V,i)+V: i=1, \ldots , n(V) \} \supseteq A$. Define 
\begin{align*}
B(V) = \{ & x^* \in X^* : \text{ there exists } 1 \leq i \leq n(V) \text{ such that } \\
& \text{Re}\, x^*x \geq 0 \text{ for all }x \in A\cap (a(V,i) + V), 
\text{ and }\text{Re}\, x^* a(V,i) =1 \} . 
\end{align*}
Since $A\cap (a(V,i) + V)$ is convex, has non-empty norm-interior and does not contain 
the origin, by the Eidelheit separation theorem, $B(V)$ is non-empty.
Obviously $B(V)$ is weakly closed. In order to see that $B(V)$ is bounded,
notice that since $a(V,i)$ belongs to the norm-interior of 
$A \cap (a(V,i) + V)$ for all $1 \leq i \leq n(V)$,
there exists $\varepsilon >0$ such that $a(V,i)-x \in A \cap(a(V,i)+V)$ for all $1 \leq i \leq n(V)$
whenever $x \in X$ with $\| x \| \leq \varepsilon$. Now let $x^* \in B(V)$. There exists 
$1 \leq i \leq n(V)$ such that $\text{Re}\, x^*x \geq 0$ for all $x \in A\cap (a(V,i) + V)$ 
and $\text{Re}\, x^* a(V,i) =1$. Hence for $x \in X$ with $\| x \| \leq \varepsilon$,
$$
0 \leq \text{Re}\, x^* (a(V,i) -x) = \text{Re}\, x^* a(V,i) - \text{Re}\, x^* x=
 1- \text{Re}\, x^* x.
$$
Thus $\text{Re}\, x^* x \leq 1$ for all $x \in X$ with $\| x \| \leq \varepsilon$, 
which proves that $B(V)$ is bounded. 
The proof continues by separating two cases:

\medskip
\noindent
{\bf Case 1:} There exists a symmetric convex weakly open neighborhood $V$ of the origin
in $X$ and sequences
$(x_n)_{n \in \N} \subseteq A$ and $(x_n^*)_{n \in \N} \subseteq B(V)$ satisfying
\begin{equation} \label{1'}
|x_m^* T^{k}x_{\ell} | < \frac{1}{m} \text{ for all integers }1 \leq k \leq \ell \leq m.
\end{equation}

Then, since $A$, $B(V)$ are bounded, weakly closed  and do not contain the origin, 
the sequences $(x_n)_n$ and $(x_n^*)_n$ are compatible and the proof of 
Remark~\ref{compatible} gives that there exists $v \in A$ such that either $Tv=0$ or 
$\overline{\text{Span}}\, \text{Orb}_T(Tv)$ is a non-trivial invariant subspace of $T$.

\medskip
\noindent
{\bf Case  2:} For every weakly open symmetric convex neighborhood $V$ of the origin of $X$
and for every sequence $(x_n)_{n \in \N} \in A^{\N}$ there exists $m \in \N$ such that
$$
B(V) \cap \bigcap_{\ell =1}^m \bigcap_{k=1}^{\ell}
\{ x^* \in X^*: | x^*T^k x_{\ell}   | < \frac{1}{m} \}
= \emptyset .
$$

Of course the last relationship implies that
\begin{equation} \label{4'}
B(V) \subseteq \bigcup_{\ell =1}^m \bigcup_{k=1}^{\ell}
\{ x^* \in X^* : |x^* T^k x_{\ell} | \geq \frac{1}{m} \} .
\end{equation}
For any $(x_n)_n \in A^{\N}$ let $Z((x_n)_n)$ denote the closed linear span of the set
$\{  T^k x_{\ell} : \text{ for } \ell \in \N \text{ and } 1 \leq k \leq  \ell \}$.

\medskip
\noindent
{\bf Claim 5:} For every $(x_n)_n \in A^{\N}$, the set $Z((x_n)_n)$ is dense in $A$.

Claim~$5$ implies  the third alternative 
of Theorem~\ref{Th:alternatives} and finishes the proof.
In order to prove Claim~$5$, note that since for every $(x_n)_n \in A^{\N}$, the set 
$Z((x_n)_n) \cap A$ is convex, it is equivalent to prove that $Z((x_n)_n)$ is weakly dense in $A$.
Thus is enough to prove that for every $y \in A$, weakly open 
neighborhood $V$ of the origin of $X$, and $(x_n)_n \in A^{\N}$, 
$Z((x_n)_n) \cap A^o \cap (y+V) \not = \emptyset$. To see this,
notice that if $Z((x_n)_n)\cap (A^o \cap (y+V)) = \emptyset$, then, since $Z((x_n)_n)$ and 
$A^o \cap (y+V)$ are convex sets and $A^o \cap (y+V)$ is a non-empty norm-open set 
(since $A^o$ is dense in $A$), by the Eidelheit separation theorem
there exists $x^* \in X^*$ such that the restriction of $x^*$ on 
$Z((x_n)_n)$ is equal to zero, $\text{Re}\, x^* x \geq 0$ for all 
$x \in A \cap (y+V)$, and $\text{Re}\, x^* x > 0$ for all 
$x \in A^o \cap (y+V)$. Let $V'$ be a weakly open convex symmetric neighborhood of the 
origin of $X$ such that $V'+V' \subseteq V$. Since 
$A \subseteq \cup \{ a(V',i)+V': 1 \leq i \leq n(V') \}$, there exists $1 \leq i \leq n(V')$ such
that $y \in a(V',i)+V'$. Thus $a(V',i) \in y + V'$, hence  
$A \cap (a(V',i)+V') \subseteq A \cap (y+V'+V') \subseteq A \cap (y+ V)$.  Therefore
$\text{Re}\, x^* x \geq 0$ for all $x \in A \cap (a(V',i)+V')$. Since 
$a(V',i) \in A^o \cap (a(V',i) + V') \subseteq A^o \cap (y+ V)$, obtain that 
$\text{Re}\, x^* a(V',i) > 0$.
Let $y^* = x^*/\text{Re}\, x^* a(V',i)$. Then $y^* \in B(V')$.
Hence (\ref{4'}) implies that there exists $m \in \N$, $\ell \in \{ 1, \ldots , m \}$ and 
$k \in \{ 1, \ldots \ell \}$ such that $| y^* T^k x_{\ell}  | \geq \frac{1}{m}$. 
In particular, $y^* T^k x_{\ell} \not = 0$ which contradicts the fact that 
the restriction of $y^*$ on  $Z((x_n)_n)$ is equal to zero,
and finishes the proof of Claim~$5$.
\end{proof} 

If $X$ is reflexive, $T \in {\mathcal L}\, (X)$ and the third alternative of the statement of 
Theorem~\ref{Th:alternatives} is valid for every closed ball $A$ of $X$, does 
$T$ have a  non-trivial invariant subspace?

\begin{Rmk} \label{genaralA}
The proof of Theorem~\ref{Th:alternatives} works if $A$ is a more general set than a 
closed ball. In fact $A$ can be taken to be any weakly compact convex subset of $X$ with 
non-empty norm interior. For such set $A$ it is not hard to see that $A^o$ is dense in $A$,
which was used in the proof.
\end{Rmk}

Immediate consequence of Theorem~\ref{Th:alternatives} is the following:

\begin{Cor}
Let $X$ be a reflexive Banach space, $T \in {\mathcal L} (X)$, and $A$ be any 
closed ball of $X$.  
Then either there exists $v \in A$ such that  $Tv=0$, or there exists $v \in A$ such that
$\overline{\text{Span}}\, \text{Orb}_T (Tv)$
is a non-trivial invariant subspace of $T$, or 
$A \subseteq \overline{\text{Span}}\, \{ T^k v: k \in \N \} $ for every $v \in A$.
\end{Cor}

\vspace{.2in}
\scriptsize{
\noindent
Department of Mathematics, University of South Carolina, Columbia, SC
29208.
giorgis@math.sc.edu

\end{document}